\documentclass{amsart} 

\usepackage{amsmath,amsthm}     
\usepackage{amssymb}            
\usepackage{euscript}           
\usepackage{enumerate,calc}
\usepackage[matrix,arrow,curve,frame]{xy}    

\xymatrixcolsep{1.9pc}                          
\xymatrixrowsep{1.9pc}
\newdir{ >}{{}*!/-5pt/\dir{>}}                  

\newtheorem{thm}[subsection]{Theorem}

\newtheorem{prop}[subsection]{Proposition}

\theoremstyle{definition}  

  {\end{list}}
  
\newcommand{\dfn}{\textbf} 

\newcommand{\mdfn}[1]{\dfn{\mathversion{bold}#1}} 

\newcommand{\tens}              {\otimes}               
\newcommand{\iso}               {\cong}  

\newcommand{\cat}{\EuScript}    

\newcommand{\cE}{{\cat E}}
\newcommand{\cF}{{\cat F}}

\newcommand{\cO}{{\cat O}}

\newcommand{\field}[1]  {\mathbb #1} 
\newcommand{\A}         {\field A}
\newcommand{\R}         {\field R}

\newcommand{\Z}         {\field Z}
\newcommand{\C}         {\field C}
\newcommand{\Q}         {\field Q}

\renewcommand{\P}         {\field P}

\DeclareMathOperator{\CH}{CH}

\DeclareMathOperator{\chara}{char}

\newcommand{\ra}{\rightarrow}                   
\newcommand{\la}{\leftarrow}                    
\newcommand{\lla}{\longleftarrow}               
\newcommand{\llla}[1]{\stackrel{#1}{\lla}}      

\newcommand{\inc}{\hookrightarrow}              

\newcommand{\blank}{-}                          

\newcommand{\rea}[1]{|{#1}|}             
\newcommand{\map}{\rightarrow}

\newcommand{\ceck}[1]{\Cech(#1)}         
\newcommand{\oceck}[1]{\Cech^{o}(#1)}    
\newcommand{\oreal}[1]{\rea{\oceck{U}}}  
\newcommand{\creal}[1]{\rea{\ceck{U}}}   

\newcommand{\Cech}{\check{C}}

\newcommand{\RP}{\R{\text{\sl P}}}
\newcommand{\CP}{\C{\text{\sl P}}}

\newcommand{\GR}{Gr}

\numberwithin{equation}{subsection}

\begin{document}

\title{Algebraic $K$-theory and sums-of-squares formulas}
\author{Daniel Dugger}
\author{Daniel C. Isaksen}

\address{Department of Mathematics\\University of Oregon\\Eugene, OR 97403}
\address{Department of Mathematics\\Wayne State University\\Detroit, MI 48202}

\email{ddugger@math.uoregon.edu}
\email{isaksen@math.wayne.edu}

\begin{abstract}
We prove a result about the existence of certain `sums-of-squares'
formulas over a field $F$.  A classical theorem uses
topological $K$-theory to prove that if such a formula exists over $\R$,
then certain powers of $2$ must divide certain binomial coefficients.
While it has been known that this result works over all characteristic
$0$ fields, the characteristic $p$ case has remained open.  In this
paper we prove the result for all fields, using algebraic $K$-theory
in place of topological $K$-theory.
\end{abstract}

\maketitle

\section{Introduction}

Let $F$ be a field.  A classical problem asks for what values of $r$,
$s$, and $n$ do there exist identities of the form
\[
\left( x_1^2 + \cdots x_r^2 \right) 
\left( y_1^2 + \cdots + y_s^2 \right) =
z_1^2 + \cdots + z_n^2
\]
in the polynomial ring $F[x_1,\ldots,x_r,y_1,\ldots,y_s]$, where the
$z_i$'s are bilinear expressions in the $x$'s and $y$'s.  Such an
identity is called a \dfn{sums-of-squares formula of type
\mdfn{$[r,s,n]$}}.  For the history of this problem, see the
expository papers \cite{L,Sh}.

The main theorem of this paper is the following:

\begin{thm}
\label{th:main}
Suppose that the characteristic of $F$ is not $2$.
If a sums-of-squares formula of type $[r,s,n]$ exists over $F$, then
$2^{\lfloor \frac{s-1}{2} \rfloor -i+1}$ divides $\binom{n}{i}$ for 
$n-r < i \leq \lfloor \frac{s-1}{2} \rfloor$.  
\end{thm}

As one specific application, the theorem shows that a formula of type
$[13,13,16]$ cannot exist over any field of characteristic not equal
to $2$. Previously this had only been known in characteristic zero.
(Note that the case $\chara(F)=2$, which is not covered by the
theorem, is rather trivial: formulas of type $[r,s,1]$ always exist).

In the case $F=\R$, the above theorem was essentially proven by Atiyah
\cite{A} as an early application of complex $K$-theory; the relevance
of Atiyah's paper to the sums-of-squares problem was only later
pointed out by Yuzvinsky \cite{Y}.  The result for characteristic zero
fields can be deduced from the case $F=\R$ by an algebraic argument
due to K.~Y. Lam and T.~Y. Lam (see \cite{Sh}).  Thus, our
contribution is the extension to fields of non-zero characteristic.
In this sense the present paper is a natural sequel to \cite{DI},
which extended another classical condition about sums-of-squares.

Our proof of Theorem~\ref{th:main}, given in Section \ref{se:main}, is
a modification of Atiyah's original argument.  The existence of a
sums-of-squares formula allows one to make conclusions about the
geometric dimension of certain algebraic vector bundles.  A
computation of algebraic $K$-theory (in fact just algebraic $K^0$),
given in Section \ref{se:K}, determines restrictions on what that
geometric dimension can be---and this yields the theorem.

Atiyah's result for $F=\R$ is actually slightly better than our
Theorem~\ref{th:main}.  The use of topological $KO$-theory rather than
complex $K$-theory yields an extra power of $2$ dividing some of the
binomial coefficients.  It seems likely that this stronger result holds in
non-zero characteristic as well and that it can be proved with
Hermitian algebraic $K$-theory.  

\subsection{Restatement of the main theorem}
\label{se:binom}
The condition on binomial coefficients from Theorem~\ref{th:main} can
be reformulated in a slightly different way.  This second formulation
surfaces often, and it's what arises naturally in our proof.  We
record it here for the reader's convenience.
Each of the following observations is a consequence of the
previous one:
\begin{itemize}
\item By repeated use of Pascal's identity
$\binom{c}{d} = \binom{c-1}{d-1} + \binom{c-1}{d}$, the number
$\binom{n+i-1}{k+i}$ equals a $\Z$-linear combination of the
numbers $\binom{n}{k+1}, \binom{n}{k+2}, \ldots, \binom{n}{k+i}$.
Similarly, $\binom{n}{k+i}$ is a $\Z$-linear
combination of 
$\binom{n}{k+1}, \binom{n+1}{k+2}, \ldots, \binom{n+i-1}{k+i}$.
\item An integer $b$ divides the numbers $\binom{n}{k+1},
\binom{n}{k+2}, \ldots, \binom{n}{k+i}$ if and only if it divides the
numbers
$\binom{n}{k+1}, \binom{n+1}{k+2}, \ldots, \binom{n+i-1}{k+i}$.
\item The series of statements 
\[ 2^N \mid \tbinom{n}{k+1}, 2^{N-1}\mid
\tbinom{n}{k+2},\ldots,2^{N-i+1} \mid \tbinom{n}{k+i} \]
is equivalent to the series of statements
\[ 2^N \mid \tbinom{n}{k+1}, 2^{N-1}\mid
\tbinom{n+1}{k+2},\ldots,2^{N-i+1} \mid \tbinom{n+i-1}{k+i}. \]
\item If $N$ is a fixed integer, then
$2^{N-i+1}$ divides $\binom{n}{i}$ 
for $n-r < i \leq N$ if and only if $2^{N-i+1}$ divides 
$\binom{r+i-1}{i}$ for $n-r < i \leq N$.
\end{itemize}

The last observation shows that Theorem
\ref{th:main} is equivalent to the theorem below.  This is the form in
which we'll actually prove the result.

\begin{thm}
\label{th:main2}
Suppose $\chara(F)\neq 2$.
If a sums-of-squares formula of type $[r,s,n]$ exists over $F$, then
$2^{\lfloor \frac{s-1}{2} \rfloor -i+1}$ divides $\binom{r+i-1}{i}$ for 
$n-r < i \leq \lfloor \frac{s-1}{2} \rfloor$.  
\end{thm}


\section{The main proof}
\label{se:main}

Let $q$ be the quadratic form on $\A^k$ defined by $q(x) = \sum_{i=1}^k x_i^2$.
A sums-of-squares formula
of type $[r,s,n]$ gives a bilinear map $\phi\colon \A^r \times \A^s \map
\A^n$ such that $q(x) q(y) = q(\phi(x,y))$.
We claim that $\phi$ induces a map
\[
f\colon \P^{s-1}-V_q \map \GR_r(\A^n)
\]
where $V_q \inc \P^{s-1}$ is the subvariety defined by $q(x)=0$ and
$\GR_r(\A^n)$ is the Grassmannian variety of $r$-planes in affine
space $\A^n$.  Given $y$ in $\P^{s-1} - V_q$, $f(y)$ is the $r$-plane
spanned by the vectors $\phi(e_1, y)$, $\phi(e_2, y)$, \ldots,
$\phi(e_r, y)$ of $\A^n$, where $e_1, \ldots, e_r$ is the standard
basis of $\A^r$.  To see that these vectors are linearly independent,
note that 
the sums-of-squares identity implies that
\[ \langle \phi(x,y),\phi(x',y) \rangle = q(y) \langle x,x' \rangle \]
for any $x$ and $x'$ in $\A^r$, where $\langle \blank,\blank \rangle$
denotes the `dot product' on $\A^k$ (for any $k$).  If one had
$\phi(x,y)=0$ then the above formula shows $q(y)\langle x,x'
\rangle=0$ for every $x'$; but since $q(y)\neq 0$, this can only
happen when $x=0$.

To verify that our description really gives a map of schemes, one can
restrict to a standard open subvariety of $\P^{s-1}$ intersected
with $\P^{s-1} - V_q$, i.e., the subvariety $U_i$ of those $y = [y_1: y_2:
\cdots: y_r]$ in $\P^{s-1}$ such that $q(y) \neq 0$ and $y_i \neq 0$.
It is clear that one has an algebraic map $U_i \ra
(\A^n)\times\cdots\times(\A^n)$ ($r$ factors) sending $y$ to
$(\phi(e_1,y), \phi(e_2,y), \ldots, \phi(e_r,y))$.  
To check that the image lies inside the Zariski
open of independent vectors (the Stiefel variety), it is enough to
verify this for the $E$-rational points over any field extension
$F\inc E$; but the argument in the previous paragraph does this.  The
Stiefel variety of course projects to the Grassmannian, and one
readily checks that the maps on the different $U_i$'s patch together.

The map $f$ has a special property we can state in terms of bundles.
Let $\eta_r$ denote the tautological $r$-plane bundle over the
Grassmannian, and let $\xi$ be the restriction to $\P^{s-1}-V_q$ of
the tautological line bundle $\cO(-1)$ of $\P^{s-1}$.  We claim the
map $\phi$ induces a map of bundles $\tilde{f}\colon r\xi \map \eta_r$
covering the map $f$.  The map $\tilde{f}$ takes an element of $r\xi$,
i.e., an element $y$ of $\P^{s-1} - V_q$ together with $r$ non-zero
scalar multiples $\alpha_1 y, \ldots, \alpha_r y$ of $y$, to the
vector $\phi(\alpha_1 e_1 + \cdots + \alpha_r e_r, y)$.  This vector
lies in the fiber of $\eta_r$ over the $r$-plane $f(y)$ of
$\GR_r(\A^n)$.

As before, to construct $\tilde{f}$ as a map of schemes one first
constructs it over the open subvarieties $U_i$ described above.  The
bundle $r\xi$ is trivial over each $U_i$, so this is straightforward.

One consequence of the existence of $\tilde{f}$ is that if we pull back
the bundle $\eta_r$ along $f$ to obtain a bundle $f^* \eta_r$ on
$\P^{s-1} - V_q$, this bundle is isomorphic to $r\xi$.

\begin{prop}
\label{pr:bundle}
If a sums-of-squares identity of type $[r,s,n]$ exists over $F$, then
there is an algebraic vector bundle $\zeta$ on $\P^{s-1}-V_q$ of rank $n-r$
such that
\[ r[\xi]+[\zeta]=n \]
as elements of 
the Grothendieck group
$K^0(\P^{s-1}-V_q)$ of locally free coherent sheaves on $\P^{s-1}-V_q$.
\end{prop}

\begin{proof}
The bundle $\eta_r$ is a subbundle of the rank $n$ trivial bundle,
which we denote by $n$.
Consider the quotient $n/\eta_r$, and set
$\zeta=f^*(n/\eta_r)$.  Since $n = [\eta_r] +
[n/\eta_r]$ in $K^0(\GR_r(\A^n))$, application of $f^*$ gives
$n = [f^* \eta_r] + [\zeta]$ in $K^0(\P^{s-1}-V_q)$.  Now recall that
$f^* \eta_r \iso r\xi$.
\end{proof}

The next task is to compute the Grothendieck group
$K^0(\P^{s-1}-V_q)$.  This becomes significantly easier if we assume
that $F$ contains a square root of $-1$.  The reason for this is made
clear in the next section.

\begin{prop}
\label{pr:K(P-V)}
Suppose that $F$ contains a square root of $-1$, and $\chara(F)\neq 2$. 
Let $c = \lfloor \frac{s-1}{2} \rfloor$.  Then
$K^0(\P^{s-1}-V_q)$ is isomorphic to 
$\Z[\nu]/(2^{c}\nu, \nu^2=-2\nu)$,
where $\nu = [\xi] - 1$ generates the reduced Grothendieck group
$\tilde{K}^0(\P^{s-1}-V_q) \iso \Z/2^{c}$.
\end{prop}

The proof of the above result will be deferred until the next section.
Note that $K^0(\P^{s-1}-V_q)$ has the same form as the complex
$K$-theory of real projective space $\RP^{s-1}$ \cite[Thm. 7.3]{Ad}.
To complete the analogy, we point out that when $F=\C$ the space
$\CP^{s-1}-V_q(\C)$ is actually homotopy equivalent to $\RP^{s-1}$
\cite[6.3]{Lw}.

By accepting the above proposition for the moment, we can finish the

\begin{proof}[Proof of Theorem~\ref{th:main2}]
Recall that one has operations $\gamma^i$ on $\tilde K^0(X)$ for any
scheme $X$ \cite[Exp.~V]{SGA6} (see also \cite{AT} for a very clear
explanation).  If $\gamma_t=1+\gamma^1 t +\gamma^2 t^2+\cdots$ denotes
the generating function, then the basic properties are:
\begin{enumerate}[(i)]
\item $\gamma_t(ab)=\gamma_t(a)\gamma_t(b)$.
\item For a line bundle $L$ on $X$ one has
$\gamma_t([L]-1)=1+t([L]-1)$. 
\item If $E$ is an algebraic vector bundle on $X$ of rank $k$ then
$\gamma^i([E]-k)=0$ for $i>k$.
\end{enumerate}
The third property follows from the preceding two via the
splitting principle.  

If a sums-of-squares identity of type $[r,s,n]$ exists over a field
$F$, then it also exists over any field containing $F$.  So we
may assume $F$ contains a square root of $-1$.  If we write
$X=\P^{s-1}-V_q$, then by Proposition~\ref{pr:bundle} there is a rank
$n-r$ bundle $\zeta$ on $X$ such that $r[\xi] + [\zeta] = n$ in $K^0(X)$.  
This may also be written as $r([\xi]-1)+([\zeta]-(n-r))=0$ 
in $\tilde K^0(X)$.  
Setting $\nu=[\xi]-1$ and applying the operation $\gamma_t$ we have
\[ \gamma_t(\nu)^r \cdot \gamma_t([\zeta]-(n-r))=1 \]
or
\[ \gamma_t([\zeta]-(n-r))=\gamma_t(\nu)^{-r}=(1+t\nu)^{-r}. 
\]
The coefficient of $t^i$ on the right-hand-side is $(-1)^i
\binom{r+i-1}{i}\nu^i$, which is the same as $-2^{i-1}\binom{r+i-1}{i}
\nu$ using the relation $\nu^2=-2\nu$.  Finally, since $\zeta$ has
rank $n-r$ we know that $\gamma^i([\zeta]-(n-r))=0$ for $i>n-r$.  In
light of Proposition~\ref{pr:K(P-V)}, this means that $2^{c}$ divides
$2^{i-1}\tbinom{r+i-1}{i}$ for $i > n-r$, where $c = \lfloor
\frac{s-1}{2} \rfloor$.  When $i-1 < c$, we can rearrange the powers
of $2$ to conclude that $2^{c-i+1}$ divides $\tbinom{r+i-1}{i}$ for
$n-r < i \leq c$.
\end{proof}


\section{$K$-theory of deleted quadrics}
\label{se:K}

The rest of the paper deals with the $K$-theoretic computation stated
in Proposition \ref{pr:K(P-V)}.

Let $Q_{n-1}\inc \P^{n}$ be the split quadric defined by one of the equations
\[ a_1b_1+\cdots+a_kb_k=0 \ (n=2k-1) \qquad\text{or}\qquad
a_1b_1+\cdots+a_kb_k+c^2=0 \ (n=2k).
\] 
Beware that in general $Q_{n-1}$ is not the same as the variety $V_q$
of the previous section.  However, if $F$ contains a square root $i$
of $-1$, then one can write $x^2+y^2=(x+iy)(x-iy)$.  After a change of
variables the quadric $V_q$ becomes isomorphic to $Q_{n-1}$.  These
`split' quadrics $Q_{n-1}$ are simpler to compute with, and we can
analyze the $K$-theory of these varieties even if $F$ does not contain
a square root of $-1$.

Write $DQ_n=\P^n-Q_{n-1}$, and let $\xi$ be the restriction to $DQ_n$ of 
the tautological line bundle $\cO(-1)$ of $\P^n$.
In this section we calculate $K^0(DQ_n)$ over any ground
field $F$ of characteristic not $2$.
Proposition~\ref{pr:K(P-V)} is an immediate corollary of this
more general result:

\begin{thm}
\label{th:K(DQ)}
Let $F$ be a field of characteristic not $2$.
The ring $K^0(DQ_n)$ is isomorphic to 
$\Z[\nu]/(2^{c}\nu, \nu^2=-2\nu)$,
where $\nu = [\xi] - 1$ generates the reduced group
$\tilde{K}^0(DQ_n) \iso \Z/2^{c}$ and $c= \lfloor \frac{n}{2} \rfloor$.
\end{thm}

Note to the reader: Below we will write $K^i(X)$ for what is usually
(but unfortunately) denoted $K_{-i}(X)$ in the algebraic $K$-theory
literature.

\subsection{Basic facts about $K$-theory}
\label{se:basicK}
Let $X$ be a scheme.  As usual $K^0(X)$ denotes the Grothendieck group
of locally free coherent sheaves, and $G_0(X)$ (also called $K'_0(X)$)
is the Grothendieck group of coherent sheaves \cite[Section~7]{Q}.
Topologically, $K^0(X)$ corresponds to the usual complex $K$-theory
functor $KU^0(\blank)$, whereas $G_0$ is something like a Borel-Moore
version of $KU$-{\it homology\/}.

Note that there is an obvious map $\alpha\colon K^0(X) \ra G_0(X)$
coming from the inclusion of locally free coherent sheaves into all
coherent sheaves.  When $X$ is nonsingular,
$\alpha$ is an isomorphism whose inverse
$\beta\colon G_0(X) \ra K^0(X)$ is constructed in the following way 
\cite[Exercise~III.6.9]{H}.  
If $\cF$ is a
coherent sheaf on $X$, there exists a resolution
\[ 0 \ra \cE_n \ra \cdots \ra \cE_0 \ra \cF \ra 0
\]
in which the $\cE_i$'s are locally free and coherent.  One defines
$\beta(\cF)=\sum_i (-1)^i [\cE_i]$.  This does not depend on the
choice of resolution, and now $\alpha\beta$ and $\beta\alpha$ are
obviously the identities.  This is `Poincare duality' for $K$-theory.

Since we will only be dealing with smooth schemes, we are now going to
blur the distinction between $G_0$ and $K^0$.  If $\cF$ is a coherent
sheaf on $X$, we will write $[\cF]$ for the class that it represents
in $K^0(X)$, although we more literally mean $\beta([\cF])$.  As an
easy exercise, check that if $i\colon U \inc X$ is an open immersion
then the image of $[\cF]$ under $i^*\colon K^0(X) \ra K^0(U)$ is the
same as $[\cF\mid_U]$.  We will use this fact often.

If $j\colon Z\inc X$ is a smooth embedding and $i\colon X-Z \inc X$ is
the complement, there is a Gysin sequence \cite[Prop.~7.3.2]{Q}
\[ 0 \la K^0(X-Z) \llla{i^*} K^0(X) \llla{j_!} K^0(Z) \la K^{-1}(X-Z)
\la \cdots 
\] 
(because of our degree conventions, $K^n$ vanishes for smooth schemes
when $n>0$).  The map $j_!$ is known as the Gysin map.  If $\cF$ is a
coherent sheaf, then $j_!([\cF])$ equals the class of its pushforward
$j_*(\cF)$ (also known as extension by zero).  Note that the
pushforward of coherent sheaves is exact for closed immersions.

\subsection{Basic facts about $\P^n$}
If $Z$ is a degree $d$ hypersurface in $\P^n$, then the structure
sheaf $\cO_Z$ can be pushed forward to $\P^n$ along the inclusion $Z
\map \P^n$; we will still write this pushforward as $\cO_Z$.  It has a
very simple resolution of the form $0 \ra \cO(-d) \ra \cO \ra \cO_Z
\ra 0$, where $\cO$ is the trivial rank 1 bundle on $\P^n$ and
$\cO(-d)$ is the $d$-fold tensor power of the tautological line bundle
$\cO(-1)$ on $\P^n$.  So $[\cO_Z]$ equals $[\cO]-[\cO(-d)]$ in
$K^0(\P^n)$.  From now on we'll write $[\cO]=1$.

Now suppose that $Z\inc \P^n$ is a complete intersection, defined by the
regular sequence of homogeneous equations $f_1,\ldots,f_r \in
k[x_0,\ldots,x_n]$.  Let $f_i$ have degree $d_i$.  The module
$k[x_0,\ldots,x_n]/(f_1,\ldots,f_r)$ is resolved by the Koszul
complex, which gives a locally free resolution of $\cO_Z$.
It follows that
\begin{equation}
\label{eq:beta}
 [\cO_Z]=(1-[\cO(-d_1)])(1-[\cO(-d_2)]) \cdots (1-[\cO(-d_r)])
\end{equation}
in $K^0(\P^n)$.  In particular, note that
\[ [\cO_{\P^i}] = (1-[\cO(-1)])^{n-i} \]
where $\P^i \inc \P^n$ is a linear subspace
because $\P^i$ is defined by $n-i$ linear equations.

One can compute that $K^0(\P^n)\iso \Z^{n+1}$, with generators
$[\cO_{\P^0}],[\cO_{\P^1}],\ldots,[\cO_{\P^n}]$ (see
\cite[Th.~8.2.1]{Q}, as one source).  
If $t=1-[\cO(-1)]$, then the previous paragraph tells us that
$K^0(\P^n)\iso \Z[t]/(t^n)$ as rings.  Here $t^k$ corresponds to
$[\cO_{\P^{n-k}}]$.

\subsection{Computations}
Let $n=2k$.  Recall that $Q_{2k-1}$ denotes the quadric in $\P^{2k}$
defined by $a_1b_1+\cdots+a_{k}b_{k}+c^2=0$.  The Chow ring
$\CH^*(Q_{2k-1})$ consists of a copy of $\Z$ in every dimension (see
\cite[Appendix A]{DI}, for example).  The generators in dimensions $k$
through $2k-1$ are represented by subvarieties of $Q_{2k-1}$ which
correspond to linear subvarieties $\P^{k-1}, \P^{k-2},\ldots,\P^{0}$
under the embedding $Q_{2k-1}\inc \P^{2k}$.  In terms of equations,
the $\P^{k-i}$ is defined by $c=b_1=\cdots=b_k=0$ together with
$0=a_k=a_{k-1}=\cdots=a_{k-i+2}$.  The generators of the Chow ring in
degrees $0$ through $k-1$ are represented by subvarieties $Z_i\inc
\P^{2k}$ ($k\leq i \leq 2k-1$), where $Z_i$ is defined by the
equations
\[ 0=b_1=b_2=\cdots=b_{2k-1-i}, \qquad
a_1b_1+\cdots+a_{k}b_{k}+c^2=0.
\]
Note that $Z_{2k-1}=Q_{2k-1}$. 

\begin{prop}  The group
$K^0(Q_{2k-1})$ is isomorphic to $\Z^{2k}$, with generators
$[\cO_{\P^0}],\ldots,[\cO_{\P^{k-1}}]$ and
$[\cO_{Z_{k}}]$, $\ldots$, $[\cO_{Z_{2k-1}}]$.  
\end{prop}

This result is classical; see \cite[Thm. 13.1]{S} and
\cite[Rem. 2.5.3]{J}.  We include a brief proof for completeness.  The
proof is similar to the computation of $\CH^*(Q_{2k-1})$ given in
\cite{DI} (which goes back at least to Hodge and Pedoe \cite{HP}), but
the details are somewhat different.  It is worth noting that to
prove Theorem \ref{th:K(DQ)} we won't actually need to know that
$K^0(Q_{2k-1})$ is free---all we'll need is the list of generators.

\begin{proof}
The argument is by induction.  When $k=1$ we have $Q_1\iso \P^1$, and
$K^0(\P^1)$ has the desired form.  Suppose $k>1$, and write
$Q=Q_{2k-1}$ and $Z=Z_{2k-2}$.  If $*$ denotes the point
$[1,0,0,\ldots,0]$ (i.e., $a_1=1$ and all other coordinates equal to
zero), consider the projection $Z-* \ra Q_{2k-3}$ which forgets $a_1$
and $b_1$.  This is a locally trivial fiber bundle with fiber $\A^1$;
hence $K^0(Q_{2k-3}) \map K^0(Z-*)$ is an isomorphism.  The closed
inclusion $Z-* \inc Q-*$ (whose complement is $Q-Z=\A^{2k-1}$) induces
a localization sequence in $K$-theory of the form
\[ 0 \la \Z \la K^0(Q-*) \llla{j_!} 
  K^0(Z-*) \llla{\delta} K^{-1}(\A^{2k-1}) \la \cdots
\]
Note that the pullback map $K^0(Q-*) \ra K^0(Q-Z)\iso \Z$ sends
$[\cO_{(Q-*)}]$ to $[\cO_{(Q-Z)}]$, which is the generator.  

We know by the computation in \cite{DI} that $\CH^*(Q-*)$ is free of
rank $2k-1$.  The Chern character isomorphism $K^0(Q-*)\tens \Q \iso
\CH^*(Q-*)\tens \Q$ \cite[Exp.~XIV]{SGA6} then shows that the rank of
$K^0(Q-*)$ is $2k-1$.  Since we know by induction that $K^0(Z-*) \iso
K^0(Q_{2k-3})$ is free of rank $2k-2$, the boundary map $\delta$ must
be zero.  Thus, $K^0(Q-*)$ is free as well.

Chasing through the isomorphisms in 
the above argument, $K^0(Q-*)$ is generated by the
classes $[\cO_{(W-*)}]$ where $W$ ranges through the subvarieties 
\begin{equation}
\label{eq:subvarieties}
\P^1,\ldots,\P^{k-1},Z_{k},\ldots,Z_{2k-1}.
\end{equation}

Finally, one analyzes the localization sequence for the closed
inclusion $* \inc Q$:
\[
0 \la K^0(Q-*) \la K^0(Q) \llla{j_!} K^0(*) \llla{\delta}
K^{-1}(Q-*) \la \cdots
\]
We know that $K^0(Q-*)$ is a free abelian group of rank $2k-1$, and we
know that the rank of $K^0(Q)$ must be $2k$ by comparing rational
$K$-theory to the rational Chow groups
using the Chern character isomorphism for $Q$.  It follows again that
$\delta$ is zero, and $K^0(Q)$ is free of rank $2k$.  The map $j_!$
takes the generator of $K^0(*)=\Z$ to $[\cO_{\P^0}]$.  This class
together with the classes $[\cO_{W}]$ where $W$ ranges over the
subvarieties in (\ref{eq:subvarieties}) are a free basis for $K^0(Q)$.
\end{proof}

\begin{proof}[Proof of Theorem~\ref{th:K(DQ)} when $n$ is even]
Set $n=2k$.
To calculate $K^0(DQ_{2k})$ we must analyze the localization sequence
\[ 0 \la K^0(DQ_{2k}) \la K^0(\P^{2k}) \llla{j_!} K^0(Q_{2k-1}).
\]
The image of $j_!: K^0(Q_{2k-1}) \ra K^0(\P^{2k})$ 
is precisely the subgroup generated by  
$[\cO_{\P^0}],\ldots,[\cO_{\P^{k-1}}]$ and
$[\cO_{Z_{k}}]$, $\ldots$, $[\cO_{Z_{2k-1}}]$.
Since $\P^i$ is a complete intersection defined by $2k-i$ linear equations, 
formula (\ref{eq:beta}) tells us that 
$[\cO_{\P^i}] = t^{2k-i}$ for $0 \leq i \leq k-1$.

Now, $Z_{2k-1}$ is a degree $2$ hypersurface in $\P^{2k}$, and so
$[\cO_{Z_{2k-1}}]$ equals $1-[\cO(-2)]$.  Note that
\[ 1-[\cO(-2)]= 2(1-[\cO(-1)])-(1-[\cO(-1)])^2 = 2t-t^2.
\] 
In a similar way one notes that $Z_i$ is a complete intersection
defined by $2k-1-i$ linear equations and one degree $2$ equation, so
formula (\ref{eq:beta}) tells us that
\[ [\cO_{Z_i}]=(1-[\cO(-1)])^{2k-1-i} \cdot (1-[\cO(-2)]) =
t^{2k-1-i}(2t-t^2).
\] 

The calculations in the previous two paragraphs imply that the kernel
of the map $K^0(\P^{2k}) \ra K^0(DQ_{2k})$ is the ideal generated by
$2t-t^2$ and $t^{k+1}$.  This ideal is equal to the ideal generated by
$2t-t^2$ and $2^k t$, so $K^0(DQ_{2k})$ is isomorphic to $\Z[t]/(2^k
t, 2t - t^2)$.  If we substitute $\nu=[\xi]-1=-t$, we find
$\nu^2=-2\nu$.

To find $\tilde{K}^0(DQ_{2k})$, we just have to take the (additive)
quotient of $K^0(DQ_{2k})$ by the subgroup generated by $1$.  This
quotient is isomorphic to $\Z/2^k$ and is generated by $\nu$.
\end{proof}

This completes the proof of Theorem~\ref{th:K(DQ)} in the case where
$n$ is even.  The computation when $n$ is odd is very similar; we will
only briefly outline the differences.  

\begin{proof}[Sketch proof of Theorem~\ref{th:K(DQ)} when $n$ is odd]
In this case $Q_{n-1}$ is defined by the equation
$a_1b_1+\cdots+a_kb_k=0$ with $k = \frac{n+1}{2}$.  The Chow ring
$\CH^*(Q_{n-1})$ consists of $\Z$ in every dimension except for $k-1$,
which is $\Z\oplus \Z$.  The `extra' generator in this dimension is
the projective space $P$ defined by $0=b_1=b_2=\cdots=b_{k-1}=a_k$.
One finds that $K^0(Q_{n-1})$ is free of rank $n+1$, on the same
generators as before plus this `extra' generator $[\cO_{P}]$.  The map
$j_!\colon K^0(Q_{n-1}) \map K^0(\P^n)$ sends both $[\cO_P]$ and
$[\cO_{\P^{k-1}}]$ to $t^{k-1}$, so this doesn't affect the
computation of $K^0(DQ_n)$.
\end{proof}


\bibliographystyle{amsalpha}

\begin{thebibliography}{SGA6}

\bibitem[Ad]{Ad} J.~F. Adams, \emph{Vector fields on spheres},
Ann. of Math. (2) {\bf 75} (1962) 603--632.

\bibitem[A]{A} M. F. Atiyah, \emph{Immersions and embeddings of
manifolds}, Topology {\bf 1} (1962) 125--132.

\bibitem[AT]{AT}
M. F. Atiyah and D. O. Tall,
{\em Group representations, $\lambda$-rings and the $J$-homomorphism},
Topology {\bf 8} (1969) 253--297.

\bibitem[DI]{DI} D. Dugger and D. C. Isaksen, \emph{The Hopf condition
for bilinear forms over arbitrary fields}, preprint, 2003.

\bibitem[H]{H} R. Hartshorne, \emph{Algebraic geometry}, Graduate
Texts in Mathematics {\bf 52}, Springer, 1977.

\bibitem[HP]{HP} W. V. D. Hodge and D. Pedoe, 
\emph{Methods of algebraic geometry}, Vol. II,
Cambridge University Press, 1952.

\bibitem[J]{J} J.~P. Jouanolou, \emph{Quelques calculs en $K$-th\'eorie
des sch\'emas},  
Algebraic K-theory, I: Higher K-theories (Proc. Conf. Battelle Memorial Inst.,
Seattle, 1972), pp. 317--335,
Lecture Notes in Mathematics {\bf 341}, Springer, 1973.

\bibitem[L]{L} K. Y. Lam, \emph{Topological methods for studying the
composition of quadratic forms}, 
Quadratic and Hermitian Forms, (Hamilton, Ont., 1983), pp. 173--192, 
Canadian Mathematical Society Conference Proceedings {\bf 4},
Amer. Math. Soc., 1984.

\bibitem[Lw]{Lw} P. S. Landweber, \emph{Fixed point free conjugations
on complex manifolds}, Ann. Math. (2) {\bf 86} (1967) 491--502.

\bibitem[Q]{Q}
D. Quillen, {\em Higher algebraic $K$-theory I},
Algebraic $K$-theory, I: Higher $K$-theories (Proc. Conf. 
Battelle Memorial Inst., Seattle, 1972), pp. 85--147, Lecture Notes in
Mathematics {\bf 341}, Springer, 1973.

\bibitem[SGA6]{SGA6}
{\em Th\'eorie des intersections et th\'eor\`eme de Riemann-Roch},
S\'eminaire de G\'eom\'etrie Alg\'ebrique du Bois-Marie 1966--1967 (SGA 6),
by P. Berthelot, A. Grothendieck, and L. Illusie,
Lecture Notes in Mathematics {\bf 225}, Springer, 1971.

\bibitem[Sh]{Sh} D. B. Shapiro, \emph{Products of sums of squares},
Expo. Math. {\bf 2} (1984) 235--261.

\bibitem[S]{S} R. G. Swan, \emph{Vector bundles, projective modules, and
the $K$-theory of spheres}, Algebraic Topology and Algebraic
$K$-theory (Princeton, 1983), pp. 432--522, 
Annals of Mathematics Studies {\bf 113},
Princeton University Press, 1987.

\bibitem[Y]{Y} S. Yuzvinsky, \emph{Orthogonal pairings of Euclidean
spaces}, Michigan  Math. J. {\bf 28} (1981) 131--145.

\end{thebibliography}

\end{document}